\documentclass[a4paper]{easychair}
\usepackage[utf8]{inputenc}
\usepackage{booktabs}

\hypersetup{
  pdftitle={Tipi: A TPTP-based theory development environment emphasizing proof analysis},
  pdfauthor={Jesse Alama},
  pdfkeywords={automated theorem proving, proof dependencies, proof analysis, theory development},
  colorlinks=true,
  linkcolor=black,
  citecolor=black,
}

\let\phi=\varphi

\newcommand{\coloneqq}[0]{\mathrel{\mathop:}=}

\newcommand{\systemname}[1]{\textsf{#1}}
\newcommand{\tipi}[0]{\systemname{Tipi}}
\newcommand{\mizar}[0]{\systemname{Mizar}}
\newcommand{\vampire}[0]{\systemname{Vampire}}
\newcommand{\eprover}[0]{\systemname{E}}
\newcommand{\provernine}[0]{\systemname{Prover9}}
\newcommand{\paradox}[0]{\systemname{Paradox}}
\newcommand{\macefour}[0]{\systemname{Mace4}}
\newcommand{\tptptox}[0]{\systemname{TPTP2X}}
\newcommand{\tptpforx}[0]{\systemname{TPTP4X}}
\newcommand{\getsymbols}[0]{\systemname{GetSymbols}}

\newcommand{\infinox}[0]{\systemname{Infinox}}

\newcommand{\term}[1]{\textbf{#1}}
\newcommand{\proves}[2]{{#1} \vdash {#2}}
\newcommand{\notproves}[2]{{#1} \not\vdash {#2}}
\newcommand{\powerset}[1]{\wp({#1})}

\theoremstyle{remark}
\newtheorem{definition}{Definition}

\newtheorem{example}{Example}

\newcommand{\tptpproblemurl}[1]{http://www.cs.miami.edu/~tptp/cgi-bin/SeeTPTP?Category=Problems\&Domain=GRA&File={#1}.p}
\newcommand{\tptpproblemlink}[1]{\href{\tptpproblemurl{#1}}{\texttt{#1}}}
\newcommand{\tptpaxiomurl}[1]{http://www.cs.miami.edu/~tptp/cgi-bin/SeeTPTP?Category=Axioms\&File={#1}.ax}
\newcommand{\tptpaxiomlink}[1]{\href{\tptpaxiomurl{#1}}{\texttt{#1}}}

\usepackage{mathabx}
\newcommand{\usedpremisesymbol}[0]{$\bigstar$}

\title{\tipi: A TPTP-based theory development environment emphasizing proof analysis}
\titlerunning{\tipi}

\author{Jesse Alama\thanks{Supported by the ESF
research project \emph{Dialogical Foundations of Semantics} within
the ESF Eurocores program \emph{LogICCC} (funded by the Portuguese
Science Foundation, FCT LogICCC/0001/2007).  Research for this paper was partially done while a visiting
fellow at the Isaac Newton Institute for the Mathematical Sciences in
the programme `Semantics \& Syntax'.  The author wishes to thank Ed Zalta and Paul Oppenheimer for stimulating the development of \tipi{} and for being guinea pigs for it.}\\\affiliation{Center for Artificial Intelligence}\\\affiliation{New University of Lisbon}\\\affiliation{Portugal}\\\affiliation{\url{j.alama@fct.unl.pt}, \url{http://centria.di.fct.unl.pt/~alama/}}}
\authorrunning{Alama}

\begin{document}

\maketitle

\begin{abstract}
  In some theory development tasks, a problem is satisfactorily solved
  once it is shown that a theorem (conjecture) is derivable from the
  background theory (premises).  Depending on one's motivations, the
  details of the derivation of the conjecture from the premises may or
  may not be important.  In some contexts, though, one wants more from
  theory development than simply derivability of the target theorems
  from the background theory.  One may want to know which premises of
  the background theory were used in the course of a proof output by
  an automated theorem prover (when a proof is available), whether
  they are all, in suitable senses, necessary (and why), whether
  alternative proofs can be found, and so forth.  The problem, then,
  is to support proof analysis in theory development; the tool
  described in this paper, \tipi, aims to provide precisely that.
\end{abstract}

\section{Introduction}\label{sec:intro}

A characteristic feature of theorem proving problems arising in theory
development is that we often do not know which premises of our
background theory are needed for a proof until we find one.  If we are
working in a stable background theory in which the axioms are fixed,
we naturally include all premises of the background theory because it
a safe estimate (perhaps an overestimate) of what is needed to solve
the problem.  We may add lemmas on top of the background theory to
help a theorem prover find a solution or to make our theory more
comprehensible.  Since computer-assisted theorem proving is beset on
all sides by intractability, any path through a formal theory
development task is constantly threatened by limitations both
practical (time, memory, patience, willpower) and theoretical
(undecidability of first-order validity).  Finding even one solution
(proof, model, etc.) is often no small feat, so declaring victory once
the first solution is found is thus quite understandable and may be
all that is wanted.

In some theory development tasks, though, we want to learn more about
our problem beyond its solvability.  This paper announces \tipi, a
tool that helps us to go beyond mere solvability of a reasoning
problem by providing support for answering such questions as:
\begin{itemize}
\item What premises of the problem were used in the solution?
\item Do other automated reasoning systems derive the conclusion from
  the same premises?
\item Are my premises consistent?  Do they admit unintended models?
\item What premises are truly needed for the conclusion?  Can we find
  multiple sets of such premises?  Is there a a ``minimal'' theory
  that derives the conclusion?
\item Are my axioms independent of one another?
\end{itemize}
Let us loosely call the investigation of these and related questions
\emph{proof analysis}.


\tipi{} is useful for theory exploration both in the context of
discovery and in the context of justification.  In the context of
discovery, one crafts lemmas, adds or deletes axioms, changes existing
axioms, modifies the problem statement, etc., with the aim of
eventually showing that the theory is adequate for one's purposes (it
derives a certain conjecture, is satisfiable, is unsatisfiable, etc.).
In the context of discovery, the set of one's axioms is in flux, and
one needs tools to help ensure that the development is not veering too
far off course into the unexpected countersatisfiability, admitting
``nonsense'' models, being inconsistent, etc.  In the context of
justification, after the initial work is done and a solution is found,
one wants to know more about the relationship between the theory and
the conjecture than simply that the latter is derivable from the former.
What is the proof like?  Are there other proofs?  \tipi{} is designed
to facilitate answering questions such as these.

The theorem provers and model finders that make up \tipi{} include
\eprover, \vampire, \provernine, \macefour, and \paradox.  The system
is extensible; adding support for new automated reasoning systems is
straightforward because we rely only on the SZS ontology to make
judgments about theorem proving problems.

\tipi{} uses a variety of automated reasoning technology to carry out
its analysis.  It uses theorem provers and model finders and is based
on the TPTP syntax for expressing reasoning
problems~\cite{sutcliffe2009tptp} and the SZS problem status
ontology~\cite{sutcliffe2008szs} thereby can flexibly use of a variety
of automated reasoning tools that support this syntax.

Going beyond solvability and demanding more of our solutions is
obviously not a new idea.  Our interests complement those of Wos and
collaborators, who are also often interested not simply in
derivability, but in finding proofs that have certain valuable
properties, such as being optimal in various senses;
see~\cite{wos2003automated,wos2005flowering}.  \tipi{} emphasizes
proof analysis at the level of sets of premises, whereas one could be
interested in more fine grained information such as the number of
symbols employed in a proof, whether short proofs are available,
whether a theory is axiomatized by a single formula, etc.  Such
analysis tends to involve rather expert knowledge of particular
problems and low-level tweaking of proof procedures.  \tipi{} uses
automated reasoning technology essentially always in ``automatic
mode''.

The philosophical background of \tipi{} is a classic problem in the
philosophy of logic known as the proof identity problem:
\begin{quote}
  \emph{When are two proofs the same?}
\end{quote}
Standard approaches to the proof identity problem work with natural
deduction derivations or category theory.  One well-known proposal is
to identify ``proof'' with a natural deduction derivation, define a
class of conversion operations on natural deduction derivations, and
declare that two proofs are the same if one can be converted to the
other.  See~\cite{dosen2003identity} for a discussion of this approach
to the proof identity problem.  The inspiration for \tipi{} is to take
on the proof identity problem with the assistance of modern automated
reasoning tools.  From this perspective, the TPTP
library~\cite{sutcliffe2009tptp} can be seen as a useful resource on
which to carry out experiments about ``practical'' proof identity.
TPTP problems typically don't contain proofs in the usual sense of the
term, but they do contain hints of proofs in the sense that they
specify axioms and perhaps some intermediate lemmas.

One does not need to share the philosophical background (or even care
about it) to start using \tipi, which in any case was designed to
facilitate TPTP-based theory development.

\emph{Terminology} In the following we sometimes equivocate on
the term \emph{theory}, understanding it sometimes in its mathematical
sense as a set of formulas closed under logical consequence (and so is
always infinite and has infinitely many axiomatizations), and
sometimes in its practice sense, represented as a TPTP problem, which
always has finitely many axioms.  ``TPTP theory'' simply means an
arbitrary (first-order) TPTP problem.  Of course, from a TPTP theory
$T$ we can obtain a theory in the mathematical sense of the term by
simply reading the formulas of $T$ as logical formulas and closing $T$
under logical consequence.  From an arbitrary first-order theory in
the mathematical sense of the term obviously one cannot extract a
unique finite axiomatization and, worse, many theories of interest are
not even finitely axiomatizable.  Still, we may at times, for
precision, need to understand ``theory'' in its mathematical sense,
even though of course we shall always work with finite TPTP theories
(problems).

\emph{Convention} Some TPTP problems do not have a conjecture formula.
Indeed, some TPTP problems are not theorem proving problems per se but
are better understood as model finding problems (e.g., the intended SZS
status is \verb+Satisfiable+).  For expository convenience we shall
restrict ourselves to TPTP problems whose intended interpretation is
that a set of premises entails a single conclusion.

The structure of this paper is as follows.  Section~\ref{sec:syntax}
describes some simple tools provided by \tipi{} to facilitate theory
development.  Section~\ref{sec:desiderata} discusses the problem of
determining which premises are needed.  Section~\ref{sec:reproving}
discusses two algorithms, one syntactic and the other semantic, for
determining needed needed premises.  Section~\ref{sec:independence}
concentrates on independent sets of axioms.
Section~\ref{sec:models} discusses some simple model analysis tools
provided by \tipi{}.  Section~\ref{sec:experience} gives a sense of
the experience so far with using \tipi{} on real-world proof analysis
tasks.  Section~\ref{sec:availability} says where can obtain \tipi{}
and briefly discusses its implementation.
Section~\ref{sec:conclusion} concludes and suggests further directions
for proof analysis and dependencies.

\section{Syntax analysis}\label{sec:syntax}

When designing TPTP theories, one needs to be careful about the
precise language (signature) that one employs.  An all-too-familiar
problem is typos: one intends to write \verb+connected_to+ but writes
\verb+conected_to+ instead.  One quick check that can help catch this
kind of error is to look for unique occurrences of constants,
functions, or predicates.  A unique occurrence of a relation or
function symbol is a sign (though by no means a necessary or
sufficient condition) that the theory is likely to be trivially
inadequate to the intended reasoning task because it will fail to be
(un)satisfiable, or fail to derive the conjecture.  Detecting such
hapax legomena early in the theory development process can prevent
wasted time ``debugging'' TPTP theories that appear to be adequate but
which are actually flawed.

\section{Needed premises}\label{sec:desiderata}

Once it is known that a conjecture $c$ is derivable from a background
theory $T$, one may want to know about the class of proofs that
witness this derivability relation.  Depending on which automated
theorem prover (ATP) is used, there may not even be a derivation of
$c$ from $T$, but only the judgment that $c$ is a theorem of $T$.  If
one does have a derivation (e.g., a resolution proof) $d$, one can
push the investigation further:
\begin{itemize}
\item Which premises of $T$ occur in $d$?
\item Are all the premises occurring in $d$ needed?
\end{itemize}
Various notions of ``needed'' are available.  For lack of space we
cannot give a complete survey of this interesting concept;
see~\cite{alama2011dependencies} for a more thorough discussion of the
notion of ``dependency'' in the context of interactive theorem
provers.  One can distinguish whether a formula is needed for a
\emph{derivation}, or for a \emph{conclusion}.  In a Hilbert style
calculus, the sequence $d$ of formulas
\[
\langle C, A, A \rightarrow B, B \rangle
\]
is a derivation of $D$ from the premises $X = \{C,A,A \rightarrow
B\}$.  All axioms of $X$ do appear in $d$, so it is reasonable to
assert that all of $X$ is needed for $d$.  But are all premises needed
for the conclusion $B$ of $d$?  The formula $C$ is not used as the
premise of any application of a rule of inference (here, modus ponens
is the only rule of inference).  Thus one can simply delete the first
term of $d$ and obtain a derivation $d^{\prime}$ of $B$ from $X - \{ C
\}$.  In a plain resolution calculus, a derivation of the empty clause
from a set $\mathcal{C}$ of clauses can have unused premises in the
same sense as there can be unused premises of a Hilbert-style
derivation.  Still, there can be ``irrelevant'' literals in clauses of
$\mathcal{C}$ whose deletion from $\mathcal{C}$ and from a refutation
$d$ of $\mathcal{C}$ yields a more focused proof.

Intuitively, any premise that is needed for a conclusion is also
needed for any derivation of the conclusion (assuming sensible notions
of soundness and completeness of the calculi in which derivations
are carried out).  However, a premise that is needed for a derivation
need not be needed for its conclusion.  Clearly, multiple proofs of a
conclusion are often available, employing different sets of premises.

In the ATP context, we may even find that, if we keep trimming unused
premises from a theory $T$ that derives a conjecture $c$ until no more
trimming is possible (so that every premise is needed by the ATP to
derive the conjecture), there may still be proper subsets of
``minimal'' premises that suffice.  The examples below in
Section~\ref{sec:experience} illustrate this.

\section{Reproving and minimal subtheories}\label{sec:reproving}

Given an ATP $A$, a background theory $T$ and a conjecture $c$, assume
that $T$ does derive $c$ and is witnessed by an $A$-derivation $d$.
Define
\[
T_{0} \coloneqq \{ \phi \in T \colon \text{$\phi$ occurs in $d$} \}
\]
as the set of premises of $T$ occurring in $d$.  Do we need all of $T$
to derive $c$?  If $T_{0}$ is a proper subset of $T$, then the answer
is evidently ``no''.  One simple method to investigate the question of
which premises of $T$ are needed to derive $c$ is to simply repeat the
invocation of $A$ using successively weaker subtheories of $T$.  Given
$T_{k}$ and an $A$-derivation $d_{k}$ of $c$ from $T_{k}$, define
\[
T_{k+1} \coloneqq \{ \phi \in T_{k} \colon \text{$\phi$ occurs in $d_{k}$} \}
\]
We are then after the fixed point of the sequence
\[
T_{0} \supset T_{1} \supset T_{2} \supset \dots
\]
We can view this discussion as the definition of a new proof procedure:
\begin{definition}[Syntactic reproving]
  Given a background theory $T$, an ATP $A$, a conjecture $c$, use $A$
  to derive $c$ from $A$.  If this succeeds, extract the
  premises of $T$ that were used by $A$ to derive $c$; call this set
  $T^{\prime}$.  If $T^{\prime} = T$, then stop and return $d$.  If
  $T^{\prime}$ is a proper subset of $T$, then let $T \coloneqq
  T^{\prime}$ and repeat.
\end{definition}

We call this proof procedure ``syntactic'' simply because we view the
task of a proof finder as a syntactic one.  The name is not ideal
because an ATP may use manifestly semantic methods in the course of
deriving $c$ from $T$.  The definition of syntactic reproving requires
of $A$ only that we can compute from a successful search for a
derivation, which premises were used; we do not require a derivation
from $A$, though in practice various ATPs do in fact emit derivations
and from these we simply extract used premises.

If $A$ is a complete ATP, then we can find a fixed point, provided we
have unlimited resources; the existence of a fixed point follows from
the assumption that $T_{0} \vdash c$ and the fact that $T_{0}$ is
finite.\footnote{There may even be multiple $A$-minimal subtheories of
  $T$ that derive $c$, but the proof procedure under discussion will
  find only one of them.}  Of course, we do place restrictions on our
proof searches, so we often cannot determine that a proper subset of
$T$ suffices to derive $c$, even if there is such a subset.

It can happen that the syntactic reprove procedure applied with an ATP
$A$, a theory $T$, and a conjecture $c$, terminates with a subtheory
$T^{\prime}$ of $T$ even though there is a proper subset
$T^{\prime\prime}$ of $T$ that also suffices and, further, $A$ can
verify that $T^{\prime\prime}$ suffices.  The syntactic reprove
procedure does not guarantee that the solution it finds is truly
minimal.  Some other proof procedure, then, is needed.

\begin{definition}[Semantic reproving]\label{def:semantic-reprove}
  Given a background theory $T$ and ATPs $A$ and $B$, a conjecture $c$,
  use $A$ to derive $c$ from $A$.  If this succeeds, extract the
  premises of $T$ that were used by $A$ to derive $c$; call this set
  $T^{\prime}$.  Define
  \[
  T^{*} \coloneqq \{ \phi \in T^{\prime} \colon \notproves{T^{\prime} - \{ \phi \}}{c} \}
  \]
  Now use $A$ to check whether $\proves{T^{*}}{c}$.  If this succeeds,
  return $T^{*}$.
\end{definition}

The semantic reprove procedure takes two ATPs $A$ and $B$ as
parameters.  $A$ is used for checking derivability, whereas $B$ is used
to check underivability. This proof procedure is called ``semantic''
because the task of constructing $T^{*}$ is carried out in \tipi{}
using a model finder (e.g., \paradox{} or \macefour), which solves the
problem of showing that $\notproves{X}{\phi}$ by producing a model of
$X \cup \{ \neg \phi \}$. As with ``syntactic'' in ``syntactic
reprove'', the ``semantic'' in ``semantic reprove'' is not ideal
because any ATP that can decide underivability judgments would work;
whether $B$ uses syntactic or semantic methods (or a combination
thereof) to arrive at its solution is immaterial.  Indeed, in
principle, for $B$ a theorem prover could be used.  Even though
\vampire{} and \eprover{} are typically used to determine
derivability, because of the properties of their search procedures,
they can be used for determining underivability, though establishing
underivability is not necessarily their strong suit and often a model
finder can give an answer more efficiently to the problem of whether
$\proves{X}{\phi}$.

\section{Independence}\label{sec:independence}

In proof analysis a natural question is whether, in a set $X$ of
axioms, there is an axiom $\phi$ that depends on the others in the
sense that $\phi$ can be derived from $X - \{ \phi \}$.

\begin{definition}[Independent set of formulas] A set $X$ of formulas
  is \term{independent} if for every formula $\phi$ in $X$ it is not
  the case that $\proves{X - \{ \phi \}}{\phi}$.
\end{definition}
\tipi{} provides a proof procedure for testing whether a set of
formulas is independent.  The algorithm for testing this is
straightforward: given a finite set $X = \{ \phi_{1}, \dots,
\phi_{n} \}$ of axioms whose independence we need to test, test
successively
\begin{enumerate}
\item whether there is any $\phi = \phi_{k}$ in $X$ such that for some $j_{1} \neq k$, we have $\proves{\{ x_{j_{1}}
    \}}{\phi}$
\item whether there is any $\phi = \phi_{k}$ in $X$ such that for some $j_{1}, j_{2} \neq k$ we have $\proves{\{ x_{j_{1}}, x_{j_{2}} \}}{\phi}$,
\item etc., for increasingly larger $k$ (the upper bound is of course
  $n - 1$).
\end{enumerate}
On the assumption that most set of axiom that arise in practice are
not independent, \tipi{} employs a ``fail fast'' heuristic: if $X$ is
not independent, then we can likely find, for some axiom $\phi$ in
$X$, a small subset $X^{*}$ suffices to prove $\phi$.  Other
algorithms for testing independence are conceivable.  It could be that
the naive algorithm that is immediately suggested by the
definition---enumerate the axioms, checking for each one whether it is
derivable from the others---may be the best approach.  Experience
shows this is indeed an efficient algorithm if one really does have an
independent set (obviously the iterative ``fail fast'' algorithm
sketched requires $n(n - 1) = O(n^{2})$ calls to an ATP for a set of
axioms of size $n$, whereas the obvious algorithm just makes $n$
calls).  A model finder can be used to facilitate this: if $X - \{
\phi \} \cup \{ \neg \phi \}$ is satisfiable, then $\phi$ is
independent of $X$.  If one is dealing with large sets of axioms,
testing independence becomes prohibitively expensive, so one could
employ a randomized algorithm: randomly choose an axiom $\phi$ and a
proper subset $T^{\prime}$ of $T$ that does not contain $\phi$ and
test whether $T^{\prime}$ proves $\phi$.  \tipi{} implements all these
algorithms for checking independence.

A typical application of independence checking first invokes one of the
minimization algorithms described in Section~\ref{sec:reproving}.  If
there is a proper subset $T^{\prime}$ of $T$ that suffices to derive
$c$, then the independence of the full theory $T$ is probably less
interesting (and in any event requires more work to determine) than
the independence of the sharper set $T^{\prime}$.

Checking independence is related to the semantic reprove algorithm
described in Section~\ref{sec:reproving}.  If we are dealing with a
theory that has a conjecture formula, then the two notions are not
congruent, because the property of independence holds for a set of
formulas without regard to whether they are coming from a theory that
has a conjecture formula.  If we are dealing with a theory $T$ without
a conjecture whose intended SZS status is \verb+Unsatisfiable+, i.e.,
the theory should be shown to be unsatisfiable, then an axiom $\phi$
of a theory $T$ gets included in the the set $T^{*}$ of
Definition~\ref{def:semantic-reprove} of semantic reproving if $T - {
  \phi }$ is satisfiable.  $T$ is ``semantically minimal'' when no
proper subtheory of $T$ is unsatisfiable, i.e., every proper subtheory
of $T$ is satisfiable.  Independence and semantic minimality thus
coincide in the setting of theories without a conjecture formula with
intended SZS status \verb+Unsatisfiable+.




\section{Model analysis}\label{sec:models}

When developing formal theories, one's axioms, lemmas, and conjecture
are typically in flux.  One may request the assistance of an automated
reasoning system to check simply whether one's premises are
consistent.  One might go further an ask whether, if one is dealing
with a TPTP theory that has a conjecture, there the theory is
satisfiable when the conjecture is taken as simply another axiom.  An
``acid test'' for whether one is proceeding down the right path at all
is whether one's problem is countersatisfiable.

\tipi{} provides tools for facilitating this kind of analysis.  A
single command is available that can check, given a theory
\begin{itemize}
\item whether the theory without the conjecture has a model
\item whether the axioms of the theory together with the conjecture
  (if present) has a model
\item whether the axioms of the theory together with the negation of
  the conjecture has a model.
\end{itemize}
The second consistency check is useful to verify that one's whole
problem (axioms together with the conjecture, considered as just
another axiom) is sensible.  It can happen that the axioms of a
problem have very simple models, but adding the conjecture makes the
models somewhat more complicated.  If a set of axioms has a finite
model but we cannot determine reasonably quickly that the axioms
together with the conjecture have a finite
model, then we can take such results as a sign that the conjecture may
not be derivable from the axioms.  (Of course, it is possible that
that the set of axioms is finitely satisfiable but the set containing
axioms and the conjecture is finitely unsatisfiable.  One can use
tools such as \infinox{}~\cite{claessen2009automated} to complement
\tipi{} in such scenarios.)

\section{Experience}\label{sec:experience}

\tipi{} has so far been used successfully to analyze a variety of
theories occurring in diverse TPTP theory development tasks.  It has
proved quite useful for theory development tasks in computational
metaphysics~\cite{fitelson2007steps}, which was the initial impetus
for \tipi.

To get a sense of how one can apply \tipi, we now consider several
applications of \tipi{} to problems coming from the large TPTP library
of automated reasoning problems.  In these examples we use \eprover{}
as our theorem prover and \paradox{} as our model finder.

\begin{example}[\tptpproblemlink{GRA008+1}]
  A problem in a first-order theory\footnote{See \tptpaxiomlink{GRA001+0}.} about graphs has 17 premises.
  Syntactic reprove brings this down to
  12. Progressing further with semantic repoving, we find that 8 of
  the 12 are needed (in the sense that for each of them, their
  deletion, while keeping the others, leads to countersatisfiability).
  Moreover, none of the other 4 is individually needed (the conjecture
  is still derivable from the 4 theories one obtains by deleting the
  4).  It turns out that there are two minimal theories that suffice
  to derive the conjecture; see Table~\ref{tab:gra008+1-minimal}.
  \begin{table}
    \centering
    \begin{tabular}{rp{5cm}cc}
      Axiom Name & Formula & Minimum 1 & Minimum 2\\
    \toprule
    \texttt{edge\_ends\_are\_vertices} & The head and the tail of an edge are vertices & \usedpremisesymbol & \\
    \midrule
    \texttt{in\_path\_properties} & If $P$ is a path from vertices $v_{1}$ and $v_{2}$, and if $v$ is in $P$, then~(i)~$v$ is a vertex and~(ii)~there is an edge $e$ of $P$ such that $v$ is either the head or the tail of $e$. & & \usedpremisesymbol \\
    \midrule
  \texttt{on\_path\_properties} & If $P$ is a path from vertices $v_{1}$ and $v_{2}$, and if $e$ is on $P$, then $e$ is an edge and both the head and the tail of $e$ are in $P$. & & \usedpremisesymbol \\
    \bottomrule
  \end{tabular}

  \caption{Two minimal subtheories of \tptpproblemlink{GRA008+1}.}
\label{tab:gra008+1-minimal}
\end{table}
\end{example}

\begin{example}[\tptpproblemlink{PUZ001+1}]
  Pelletier's Dreadbury Mansion puzzle~\cite{Pel86-JAR} asks: ``Who killed Aunt
  Agatha?''
  \begin{quotation}
     Someone who lives in Dreadbury Mansion killed Aunt Agatha.
     Agatha, the butler, and Charles live in Dreadbury Mansion,
     and are the only people who live therein. A killer always
     hates his victim, and is never richer than his victim.
     Charles hates no one that Aunt Agatha hates. Agatha hates
     everyone except the butler. The butler hates everyone not
     richer than Aunt Agatha. The butler hates everyone Aunt
     Agatha hates. No one hates everyone. Agatha is not the
     butler.
  \end{quotation}
  Among the $13$ first-order sentences in the formalization are three
  \begin{quotation}
    \texttt{lives(agatha)}, \texttt{lives(butler)}, \texttt{lives(charles)}
  \end{quotation}
  that turn out to be deletable, which the reader may find amusing
  since we are dealing with a murder mystery.  Each of the $10$ other
  premises turn out to be needed (their deletion leads to
  countersatisfiability), so there is a unique minimal subtheory of
  the original theory that suffices to solve the mystery (which is
  that Agatha killed herself).  Note that the premise that Agatha is
  not the butler (an inference that would perhaps be licensed on
  pragmatic grounds if it were missing from the text) is needed, which
  perhaps explains why the puzzle explicitly states it.
\end{example}

\begin{example}[\tptpproblemlink{REL002+1}]
  A problem about relation algebra is to show that $\top$ is a right
  unit for the join operation ($+$):
  \[
  \forall x (x \vee \top = \top) .
  \]
  There are 13 premises\footnote{See \tptpaxiomlink{REL001+0}.}.
  Syntactic reprove with \vampire{} shows that $7$ axioms can be cut,
  whereas syntactic reprove with \eprover{} finds $6$.  The sets of
  syntactically minimal premises of \eprover{} and \vampire{} are,
  interestingly, not comparable (neither is a subset of the other).
  Semantic reprove with the 10 distinct axioms used by either
  \vampire{} or \eprover{} shows, surprisingly, that $2$ are needed
  whereas each of the other $8$ is separately eliminable.  Of the 256
  combinations of these $8$ premises, we find two minima; see
  Table~\ref{tab:rel002+1-minimal}.

  \begin{table}
    \centering
    \begin{tabular}{rp{3cm}cc}
      Axiom Name & Formula & Solution 1 & Solution 2\\
      \toprule
      \texttt{composition\_identity} & $x;\mathbf{1} = x$ & \usedpremisesymbol & \\
      \midrule
      \texttt{converse\_cancellativity} & $\breve{x};\overline{x;y} \vee \bar{y} = \bar{y}$ & \usedpremisesymbol & \\
      \midrule
      \texttt{converse\_idempotence} & $\bar{\bar{x}} = x$ & \usedpremisesymbol & \\
      \midrule
      \texttt{converse\_multiplicativity} & $ \breve{x;y}  =  \breve{x};\breve{y}  $ & \usedpremisesymbol & \\
      \midrule
      \texttt{maddux3\_a\_kind\_of\_de\_Morgan} & $x = \overline{\overline{x} \vee \overline{y}} \vee \overline{\overline{x} \vee y}$ & & \usedpremisesymbol \\
      \bottomrule
    \end{tabular}
  \caption{Minimal subtheories of \tptpproblemlink{REL002+1}.}
\label{tab:rel002+1-minimal}
\end{table}

\end{example}

%

With enough caution, \tipi{} can be used somewhat in the large-theory
context, where there are ``large'' numbers of axioms (at least several
dozen, sometimes many more). Although it is quite hopeless, in the
large-theory context, to test all possible combinations of premises in
the hope of discovering all minimal theories, one can, sometimes use
syntactic reproving to weed out large classes of subsets.  With these
filtered premises, semantic reproving can be used to find minima using
a more tractable number of combinations of premises.

\begin{example}[\tptpproblemlink{TOP024+1}]
  Urban's mapping~\cite{urban2006mptp} of the \mizar{} Mathematical
  Library, with its rich language for interactively developing
  mathematics, into pure first-order theorem proving problems, is a
  rich vein of theorem proving problems.  Many of them are quite
  challenging owing to the large number of axioms and the inherent
  difficulty of reasoning in advanced pure mathematics.

  Here the problem is to prove that every maximal $T_{0}$ subset of a
  topological space $T$ is dense.

  Of the $68$ available premises, 9 are found through an initial
  syntactic reproving run using \eprover{} and \vampire{}.  Of these
  9, 3 are (separately) not needed, whereas the other $6$ are needed.
  Of the $8$ combinations of these $3$ premises, we find two minima;
  see Table~\ref{tab:TOP024+1-minimal}.

  \begin{table}
    \centering
    \begin{tabular}{rp{6cm}cc}
      Name & Formula & Solution 1 & Solution 2\\
      \toprule
      \texttt{dt\_k3\_tex\_4} & Maximal anti-discrete subsets of a topological space $T$ are subsets of $T$ & \usedpremisesymbol & \\
      \midrule
      \texttt{reflexivity\_r1\_tarski} & $X \subseteq X$ & & \usedpremisesymbol \\
      \midrule
      \texttt{t3\_subset} & $A \in \powerset{X}$ iff $A \subseteq X$ & & \usedpremisesymbol \\
      \bottomrule
  \end{tabular}
  \caption{Two minimal subtheories of \tptpproblemlink{TOP024+1}.}
\label{tab:TOP024+1-minimal}
\end{table}

\end{example}


\section{Availability}\label{sec:availability}

\tipi{} is available at
\begin{quote}
  \url{https://github.com/jessealama/tipi}
\end{quote}
At present \tipi{} relies on the \getsymbols{}, \tptptox, and
\tptpforx{} tools, which are part of the TPTP World
distribution~\cite{sutcliffe2010tptp}.  These are used to parse TPTP
theory files; a standalone parser for the TPTP language is planned,
which would eliminate the dependency on these additional tools.

\section{Conclusion and future work}\label{sec:conclusion}

At the moment \tipi{} supports a handful of theorem provers and model
finders.  Supporting further systems is desirable; any automated
reasoning system that supports the SZS ontology could, in principle,
be added.

\tipi{} supports, at the moment, only first-order logic, and so covers
only a part of the space of all TPTP theories.  There seems to be no
inherent obstacle to extending \tipi{} to support higher-order
theories as well.

More systematic investigation for alternative proofs of a theorem
could be carried out using \provernine's clause weight mechanism.  One
could have an alternative approach to the problem of generating
multiple alternative proofs to the simple approach taken by \tipi.

When working with models of a theory under development that makes true
some rather unusual or unexpected formulas, it can sometimes be
difficult to pinpoint the difficulty with the theory that allows it to
have such unusual models.  One has to infer, by looking at the raw
presentation of the model, what the strange properties are.  We would
like to implement a smarter, more interactive diagnosis of ``broken'' theories.

The problem of finding minimal subtheories sufficient to derive a
conjecture, checking independence of sets of axioms, etc., clearly
requires much more effort than simply deriving the conjecture.
\tipi{} thus understandably can take a lot of time to answer some of
the questions put to it.  Some of this inefficiency seems unavoidable,
but it is reasonable to expect that further experience with \tipi{}
could lead to new insights into the problem of finding theory minima,
determining independence, etc.

The proof procedures defined by \tipi{} naturally suggest extensions
to the SZS ontology~\cite{sutcliffe2008szs}.  One can imagine SZS statuses such as

\begin{itemize}
\item \verb+IndependentAxioms+: The set of axioms is independent.
\item \verb+DependentAxioms+: The set of axioms is dependent.
\item \verb+MinimalPremises+: No proper subset of the axioms suffices
  to derive the conjecture.
\item \verb+NonMinimalPremises+: A proper subset of the
  axioms suffices to derive the conjecture.
\item \verb+UniqueMinimum+: There is a unique subset $S$ of the axioms such that $S$ derives the conjecture and every proper subset of $S$ fails to derive the conjecture.
\item \verb+MultipleIncomparableMinima+: There are at least two proper
  subsets $S_{1}$ and $S_{2}$ of the axioms suffices to derive the
  conjecture, with neither $S_{1} \subseteq S_{2}$ nor $S_{2} \subseteq
  S_{1}$.
\end{itemize}
\tipi{} itself already can be seen as supporting these (currently
unofficial) SZS statuses.  One could even annotate the statistics for
many problems in the TPTP library by listing the number of possible
solutions (minimal subtheories of the original theory) they admit, or
the number of premises that are actually needed.

\bibliographystyle{plain}
\bibliography{tipi}

\end{document}